\documentstyle[art10,twoside]{article}
\catcode`@=11
\def\ps@myheadings{\let\@mkboth\@gobbletwo
 \def\@oddhead{\hfil{\sl\rightmark}\hfil}%
 \def\@oddfoot{\hfil\rm\thepage}\def\@evenhead{\hfil{\sl\leftmark}\hfil}%
 \def\@evenfoot{\rm\thepage\hfil}\def\sectionmark##1{}
\def\subsectionmark##1{}}
\catcode`@=12
\newcommand{\beq}{\begin{equation}}
\newcommand{\eeq}{\end{equation}}
\newcommand{\nit}{\noindent} 
\newcommand{\vs}{}

\newcommand{\der}{{{\cal D}}}

\begin{document}

%\pagestyle{plain}
%\pagenumbering{arabic}

%\begin{flushright}
%DAMTP-96/75 \\
%\end{flushright}

\begin{center}
{\bf THE $q$-CALCULUS FOR GENERIC $q$ AND $q$ A ROOT OF UNITY} 
\footnote{Presented  
at $5^{th}$ Colloquium `Quantum groups and integrable 
systems', Prague, 20-22 June 1996.}
\\[0.3ex]

{\sc R.S. \ Dunne, A.J. \ Macfarlane, } \\[0.3ex]

{\small{\sl DAMTP, Cambridge University, Silver Street, Cambridge, U.K.}}
\\[0.3ex]

{\sl \&} 
\\[0.3ex]

{\sc  J.A. \ de Azc\'{a}rraga and J.C. \ P\'{e}rez Bueno 
\footnote{E-mails:  r.s.dunne@damtp.cam.ac.uk, a.j.macfarlane@damtp.cam.ac.uk, 
\\
azcarrag@evalvx.ific.uv.es and pbueno@lie.ific.uv.es.}} \\[0.3ex]

\small{\sl Departamento de F\'{\i}sica Te\'{o}rica and IFIC,} \\
\small{\sl Centro Mixto CSIC-Universidad de Valencia,} \\
\small{\sl E-46100-Burjassot (Valencia) Spain.} \\[0.3ex]
\end{center}
\noindent
{\small
The $q$-calculus for generic $q$ is developed and related to the deformed 
oscillator of parameter $q^{1/2}$. By passing with care to the limit in which
$q$ is a root of unity, one uncovers the full algebraic structure of 
${{\cal Z}}_n$-graded fractional supersymmetry and its natural 
representation.}

\section{Introduction.}

\vs

The $q$-calculus, in the generic case in which the fixed complex number $q$ is
not a root of unity but is otherwise arbitrary, involves a single
non-commuting variable $\theta$ 
and its left derivative operator $\der = \der_L$ and is 
governed by the commutation relation
\beq 
\der \theta -q \theta \der =1 \quad .
\eeq
It is closely related to the $q$-deformed oscillator \cite{AC} \cite{AJM}
\cite{LCB}, as is shown below. 

The context in which $q$ is a root of unity, $q=\exp  {{2\pi i} \over n}$,
is also of great interest. It involves $\theta$ such that $\theta^n=0$ and
and can be discussed by truncating the generic case so as to exclude powers of
$\theta$ higher than the $(n-1)$-th. 
However, if 
we pass with care from the generic case to the 
limit in which $q$ is a root of unity much more structure can be exposed.
The algebraic structure in question is the full algebraic structure of 
fractional supersymmetry (FSUSY), not only the generalised 
Grassmann sector of this
${{\cal Z}}_n$-graded theory which is the part that where
$\theta$ enters but also its bosonic sector. 
The paper shows how both these `sectors' emerge and 
discusses the representation of the theory in a product Hilbert space. This has
an ordinary oscillator factor for the bosonic degree of freedom, and 
relates the generalised Grassmann sector to the $q$-deformed 
oscillator with deformation parameter 
$q^{1/2}$, which is exactly what is needed to ensure proper hermiticity
properties. 
We do not here make any extensive discussion of the interplay between 
the sectors. But some idea of the 
insights regarding this interplay can be obtained from \cite{DMdAPBplb}
which is devoted to the case of $q=-1$, which is that of ordinary 
({\it i.e.}, ${{\cal Z}}_2$-graded) supersymmetry in zero space dimension. 
It seems worthwhile emphasising that the
$q$-deformed oscillators at deformation parameter $q^{1/2}$ emerge as
those generalisations from $n=2$ to higher $n$ of the fermions
of supersymmetry which are best suited to the development of FSUSY.

References to FSUSY, including many to the extensive work of others, can be
found in our published \cite{dAM} and forthcoming \cite{DMdAPB} work.

%\newpage

\section{The $q$-calculus.}
\vs

For any graded algebra, we define a graded bracket, initially for elements 
$A$ and $B$ of pure grade $g(A)$ and $g(B)$, by 
\beq
[A,B]_{\gamma (A,B)}=AB-\gamma (A,B) \, BA, \quad
\gamma (A,B):=q^{-g(A) g(B)} \quad .
\label{gb}
\eeq
\nit This satisfies
\beq
 [AB ,C]_{\gamma (AB,C)}  =  A[B,C]_{\gamma (B,C)}+\gamma (B,C)[A,C]_
{\gamma (A,C)} B \quad , \nonumber 
\eeq
\beq
 [A ,BC]_{\gamma (A,BC)}  =  [A,B]_{\gamma (A,B)} C +\gamma (A,B) B 
[A,C]_{\gamma (A,C)} \quad ,
\eeq
\nit wherein  $g(AB) =g(A)+g(B)$ is implicit. 
The definition (\ref{gb}) extends by linearity to elements 
of the algebra not of pure grade.

In (\ref{gb}) and until section four, 
$q$ is `generic' {\it i.e.}, it is a fixed but arbitrary complex 
quantity that is not a root of $1$. To define the $q$-calculus 
algebra, we employ a single non-commuting variable $\theta$
of grade $1$, together with left and right 
derivatives $\der_L$ and $\der_R$ of grade $-1$. 
Since we shall not refer to $\der_R$ here (cf. \cite{dAM,NEW}),
we shall write $\der_L \equiv \der $.
The action 
of $\der$ upon powers 
and hence functions of $\theta$ is defined algebraically with the help
of the graded bracket
\beq
1=[\der , \theta ]_q \, := \, \der \theta -q \theta \der \quad ,
\label{deriv}
\eeq
\nit so that, for any positive integer $m$, we have
\beq
[\der ,\theta^{(m)}]_{q^m}=\theta^{(m-1)} \quad ,
\eeq
\nit 
where the bracketed exponent is defined by 
\beq
B^{(m)}:=B^m/[m]_q! \quad , \quad [m]_q=(1-q^m)/(1-q) \quad .
\label{notn}
\eeq
\nit The action extends obviously to $f(\theta)=\sum_{m=0}^\infty C_m 
\theta^{(m)}$, where the $C_m$ are complex numbers 
\beq
{{df} \over {d\theta}}
 \equiv [\der ,f(\theta)]_\gamma \quad 
= \sum_{m=0}^\infty C_m [\der ,\theta^{(m)}]_\gamma \quad  = \sum_{m=1}^\infty
C_m \theta^{(m-1)} \quad .
\eeq
\nit It extends further also to $f(\theta) =\sum_{m=0}^\infty \theta^{(m)} 
A_m$, where
the $A_m$ are quantities independent of $\theta$ and of pure grade $g(A_m)$,
so that $[\theta ,A_m]_\gamma$ with $\gamma = q^{-g(A_m)}$, giving 
\begin{eqnarray}
{{df} \over {d\theta}} 
& \equiv & [\der ,f(\theta)]_\gamma \quad 
= \sum_{m=0}^\infty [\der ,\theta^{(m)} A_m]_\gamma \quad , \nonumber \\
& = & \sum_{m=0}^\infty ( [\der ,\theta^{(m)}]_{q^m} \, A_m + q^m A_m
[\der ,A_m]_\gamma )  \quad = \sum_{m=1}^\infty \theta^{(m-1)} A_m \quad , 
\end{eqnarray}
\nit provided that $[\der , A_m]_\gamma =0$
where $\gamma =q^{g(A_m)}$ (cf. \ref{gb}), which
can be seen to be compatible with the corresponding result assumed for
$\theta$ and $A_m$. An illustration, featuring $\exp _q(\theta A)
=\sum_{m=0}^\infty (\theta A)^{(m)}$ and wherein $A$ is of pure grade, 
gives rise to
\beq 
{{d\exp _q(\theta A)} \over {d\theta}}
=A \, \exp _q(\theta A) \quad .
\eeq
\nit Another notable result, involving a parameter $\varepsilon$ of grade $1$
with $[\der ,\varepsilon ]_{q^{-1}}=0$ and
$[\varepsilon , \theta ]_{q^{-1}}=0 $,
\nit shows that the quantity $G_L(\varepsilon)=\sum_{m=0}^\infty 
\varepsilon^{(m)} \der^m $ generates the translation $\theta \mapsto \theta+
\varepsilon$, {\it i.e.}
$G_L(\varepsilon) \, f(\theta) \, G_L(\varepsilon)^{-1} =f(\theta + 
\varepsilon) \quad . $

%newpage

\section{Number, creation and destruction operators.}

\vs

We begin by constructing a number operator $N$ of grade zero 
with the properties
\begin{eqnarray}
{[}N,\theta {]}=\theta \quad , & & q^N \theta q^{-N}=q\theta \quad , 
\nonumber \\
 {[}N,\der {]}=-\der \quad , & & q^N \der q^{-N}=q^{-1}\der \quad . 
\label{numb}
\end{eqnarray}
\nit Since $N$ is of grade zero, an expression for it in terms of $\theta$
and $\der$ may be expected to be of the form
\beq
N=\sum_{m=0}^\infty C_m \theta^m \der^m \quad ,
\eeq
and 
\beq
N=\sum_{m=1}^\infty {{(1-q)^m} \over {1-q^m}} \theta^m \der^m \quad ,
\eeq
\nit satisfies both lines of (\ref{numb}). Likewise the right entries of
(\ref{numb}) may be shown to be satisfied by 
\beq
q^N=\der \theta -\theta \der \quad = 1-(1-q)\theta \der \quad ,
\label{qton}
\eeq 
\nit where (\ref{deriv}) has been used. Useful consequences of these results 
include
\beq
\theta \der = [N]_q \quad , \quad 
\theta^m \der^m = {{[N]_q!} \over {[N-m]_q!}} 
\quad , \quad
N=\sum_{m=1}^\infty {{(1-q)^m} \over {1-q^m}} {{[N]_q!} \over {[N-m]_q!}} 
\quad .
\label{props}
\eeq
\nit The last result in (\ref{props}) is of interest as it gives 
an expression for $N$ in terms of $q^N$. 
Acting on an eigenstate of $N$ whose eigenvalue is a positive
integer $r$, this yields the identity
\beq
r=\sum_{m=1}^r {{(1-q)^m} \over {1-q^m}} {{[r]_q!} \over {[r-m]_q!}} \quad .
\eeq

If we now make the identification (to within a similarity transformation,
discussion of which may be sought in \cite{DMdAPB} )
\beq
\theta =a^{\dagger}, \quad \der=q^{N/2} a \quad ,
\label{repaadag}
\eeq 
\nit then (\ref{deriv}) and (\ref{qton}) imply
\beq
a a^{\dagger}-q^{\mp 1/2} a^{\dagger} a = q^{\pm N/2} \quad.
\label{defcr}
\eeq
\nit This important result indicates how the $q$-calculus is related to the 
$q$-deformed harmonic oscillator \cite{AC}, \cite{AJM} and \cite{LCB}. If $q$ 
is real, (\ref{defcr}) admits representations in which
$a^{\dagger}$ is indeed the adjoint of $a$ in a positive definite Hilbert
space. Further, for $q=\exp  {2\pi i/n}$ when $n$ is an odd integer, the
situation to be concentrated upon below, a similar statement also holds true 
because of the fact that
the deformation parameter in (\ref{defcr}) is $q^{1/2}$. 
For simplicity the remaining sections of the Colloquium
talk confined discussion to the indicated set of roots of unity.
But the case $q=-1$ can also be treated in a similar 
spirit. It is of interest because
it underlies an instructive view \cite{DMdAPBplb} of ordinary supersymmetry 
in much the same way as the present work does for fractional supersymmetry.

\vs

%\newpage 

\section{Lemmas for use at $q$ a root of 1.}

\vs

We now confine attention 
--as mentioned at the end of the previous section--
to the $q$-values $q=\exp  {{2\pi i} \over n}$ 
for odd integer $n$. We use the shorthand ${{\cal L}}$ 
to indicate the passage to the limit in which $q$ takes on such $q$-values,
${{\cal L}}:=\lim_{q\to\exp(2\pi i/n)}$. 
We here
deduce a sequence of lemmas to be used in subsequent sections to effect the
systematic passage to the limit in question in the work of previous sections.

The following results can be proved in the order given:
\beq
{{\cal L}} {{[rn]_q} \over {[n]_q}}=r \;  , \;
{\rm for} \; {\rm integer}  \;  r \quad ,  
\quad 
{{\cal L}} {{[rn]_q!} \over {[n]_q! [(r-1)n]_q!}}=r \quad , \quad 
{{\cal L}} {{[rn]_q!} \over {([n]_q!)^r}} =r! \quad ,
\eeq 

In the next section, we shall retain $\theta^{(m)}$, in the notation
(\ref{notn}), for $m=1,2, \dots ,(n-1)$ as ${{\cal Z}}_n$- graded variables,
and explain the use of the case $m=n$ to define a variable $z$ of zero grade
by setting $z={{\cal L}} \theta^{(n)}$.
To handle powers $m$ greater than $n$, we require further lemmas to be deduced
in order. Set $m=rn+p$ for integer $r$ and $p=1,2, \dots,(n-1)$. Then we have
\beq
[rn+p]_q=[p]_q \quad , \quad
{{\cal L}} {{[rn+p]_q!} \over {[rn]_q!}}=[p]_q! \quad , \quad p=1,\ldots, 
(n-1)\quad, 
\eeq
\beq
{{\cal L}} \theta^{(rn+p)}={{\cal L}} {{\theta^{rn+p}} \over {([n]_q!)^r r!}}
\Big / {{\cal L}} {{[rn+p]_q!} \over {[rn]_q!}} \, = \, {{\theta^p} \over {[p]_q!}} \,
{1 \over {r!}} \, {{\cal L}}
\Bigl( {{\theta^n} \over {[n]_q!}} \Bigr)^r \, = \,
{{z^r} \over {r!}} \theta^{(p)} 
\quad .
\label{lem6}
\eeq
An illustration of the use of these lemmas indicates what happens to 
$\exp _q (C \theta )$, $C$ a complex number,
 in the limit under study. We find
\beq
\exp _q (C \theta )=\sum_{m=0}^\infty C^m \theta^{(m)}
= \bigl( \sum_{r=0}^\infty {{(zC^n)^r} \over {r!}} \bigr)
\bigl( \sum_{p=0}^{n-1} C^p \theta^{(p)} \bigr)\quad,
\eeq
\nit In other words
\beq 
\exp _q (C \theta )=\exp  (zC^n) \times {\rm truncated} \;
{\rm series} \quad .
\eeq

\vs

\section{The $q$-calculus for $q=\exp  {{2\pi i} \over n}$ 
for odd integer $n$.}

\vs

Now we consider what happens to the $q$-calculus for those values of $q$.
We look first at an identity valid for generic complex values of $q$ and 
any positive integer $m$, namely
\beq
[\der , \theta^{(m)} ]=\theta^{(m-1)} \quad ,
\label{idy}
\eeq 
\nit where the notation (\ref{notn}) is used.
This makes sense for $m=n$ and $[n]_q=0$ only if $\theta^n=0$ at this $q$, and
if ${{\cal L}} \theta^{(n)}$ 
attains a finite non-zero value. 
Indeed, we hereby define a new variable
$z = {{\cal L}} \theta^{(n)}$ of grade zero, so that (\ref{idy}) assumes the form
\beq 
[\der , z]=\theta^{(n-1)} \quad .
\eeq
\nit Also we see that the $q$-calculus involves the variables
\beq
1,\theta , \theta^{(2)} ,\dots , \theta^{(n-1)} \quad 
{\rm of} \; {\rm  grades} \quad
0, 1, 2, \dots ,n-1 \quad .
\label{vbles}
\eeq
\nit It is natural at this point to ask what happens to powers of 
the generalised Grassmann variable $\theta$
higher than the $n$-th. If they are simply discarded much insight into the 
nature of fractional supersymmetry \cite{DMdAPB} (and likewise of ordinary
supersymmetry \cite{DMdAPBplb}) is lost. Actually lemma (\ref{lem6}) of the
previous section gives us directly an explicit non-trivial answer to 
the question. It 
follows that the generalised superfields of the context are linear 
combinations of the variables (\ref{vbles})
with coefficients that are functions of $z$
Thus $z$ plays the role for the present (${{\cal Z}}_n$-graded 
fractional supersymmetry) context that $t$ plays in ordinary 
(${{\cal Z}}_2$-graded) supersymmetric mechanics.

Next it is natural to ask about $\der^n$ and to ask how ${{\partial} \over
{\partial z}}$ enters the picture, plainly not unrelated matters. By looking at
a suitable $n$-fold graded bracket involving $\theta$ and 
$\der$ each $n$ times,
it is not hard to  show that $\der^n$ must be a well defined
quantity such that
\beq 
[\der^n , z]=1 \quad .
\eeq
\nit Thus we make the identification $\der^n={{\partial} \over
{\partial z}}$. 

It is clear that we must adjust somewhat our view of the nature of the 
derivative operator $\der$. Presenting (\ref{idy}) in the form
\beq
[\der, z]= \theta^{(n-1)}= \Bigl( {{dz} \over {d\theta}} \Bigr)\quad ,
\eeq 
\nit suggests that we now must view $\der$ as a total derivative with 
respect to $\theta$ and write 
\beq
\der = \partial_{\theta} +\theta^{(n-1)} \partial_z \quad ,
\label{supe}
\eeq
\nit which corresponds to the result
\beq
\Bigl( {{df} \over {d\theta}} \Bigr) =\Bigl( {{\partial f} \over 
{\partial \theta}} \Bigr)+ \Bigl( {{dz} \over {d\theta}} \Bigr)
{{\partial f} \over {\partial z}} \quad .
\eeq
\nit It follows from (\ref{supe}) that
\beq
1=[\partial_{\theta}, \theta]_q \quad, \quad
(\partial_{\theta})^n=0 \quad . 
\eeq

It might be judged from the form of (\ref{supe}) that $\der$ is closely related
to the full supercharge of the ${{\cal Z}}_n$-graded fractional 
supersymmetry (FSUSY), 
and it can be seen in \cite{dAM} \cite{DMdAPB} (see \cite{NEW}
for ${\cal Z}_3$) that this is
exactly correct. That $\der$ should therefore generate the full translational
invariance of the theory is one aspect of this. We wish to exhibit
how this emerges
from the results at the end of section two where $\der$ is seen to generate
translation of $\theta$ at generic $q$. First we note that $\theta \mapsto
\theta + \varepsilon$ is compatible with $\theta^n=0$ only if $\varepsilon^n
=0,$ holds in addition, of course, to $\varepsilon \theta = q^{-1} \theta
\varepsilon$. Next, using lemmas from section four, we deduce
\beq
G_L = {{\cal L}} \, \sum_{m=0}^\infty \varepsilon^{(m)} \der^m 
= {{\cal L}} \, \sum_{r=0}^\infty \sum_{p=0}^{n-1} 
{{\varepsilon^{p} \der^p} \over {[p]_q!}} \times
{{\varepsilon^{rn} \der^{rn}} \over {([n]_q!)^r r!}} \quad .
\label{gen}
\eeq
\nit This makes it clear that we should define a grade zero parameter 
to associate with a 
translation of $z$ by means of
\beq
{{\cal L}} \varepsilon^{(n)}= z_{\varepsilon} \quad. 
\eeq
\nit For then it follows that we may write (\ref{gen}) in the form
\beq
G_L(z_{\varepsilon}, \epsilon )= \sum_{r=0}^\infty \sum_{p=0}^{n-1}
{{z_{\varepsilon}^r \partial_z^r} \over {r!}}  \varepsilon^{(p)} \der^p
=\exp (z_{\varepsilon} \partial_z) \sum_{p=0}^{n-1} 
\varepsilon^{(p)} \der^p \quad .
\eeq
\nit The first factor --an ordinary exponential of zero grade quantities--
generates $z \mapsto z+z_\varepsilon$ and the second factor is exactly
the one obtained in \cite{dAM} as the generator of translations of $\theta$ 
in the FSUSY context. However, the key result, showing that the full 
non-trivial FSUSY transformation of $z$ is generated by 
$G_L(z_{\varepsilon}, \epsilon )$, is
\beq
z \mapsto G_L z {G_L}^{-1} = 
z+z_{\varepsilon}+ \sum_{p=1}^{n-1} \varepsilon^{(p)} \theta^{(n-p)}
\quad ,
\eeq
in agreement with \cite{dAM}.
 
%newpage

\section{Reduction of the Representation space.}

\vs

It is rather obvious how we are to represent the algebra of $z, \partial_z, 
\theta$ and $\partial_{\theta}$. The first two describe a bosonic degree of 
freedom that commutes with the latter pair, 
one that describes in  Bargmann
style a harmonic oscillator Hilbert space ${{\cal V}}_{HO}$. 
Also, with the evident
analogue
\beq
\theta=a^{\dagger} \quad , \quad \partial_{\theta} =q^{N/2} \, a \quad ,
\eeq
\nit of (\ref{repaadag}), we see that (\ref{defcr}) still follows.
Also
$\theta^n=0$ and $\partial_{\theta}^n=0$ imply $a^n=0 \, , \, a^{\dagger n}=0$
so that the variables of non-zero grade are represented in a vector space
${{\cal V}}^n$ of $n$ degrees of freedom. Crucially, since (\ref{defcr}) 
involves the deformation parameter $q^{1/2}$, 
in the natural representation of $a$ and
$a^{\dagger}$ in ${{\cal V}}^n$ of positive definite metric,
the latter operator is indeed the true adjoint of the former. 

It is our purpose now to demonstrate how the structure just described emerges
from the work of section three when one passes from the case of generic $q$ to
$q=\exp  (2\pi i/n)$ for odd integer $n$. A representation of
$\der$ and $\theta$ 
at generic $q$ in a space spanned by eigenkets of $N$, namely $|m\rangle$ for
$m=0,1,2, \dots $, can be taken to within equivalence in the form
\beq
\der |m\rangle =|m-1\rangle \quad , \quad  \der |0\rangle =0 \quad ,
\quad \theta |m\rangle =[m+1]_q|m+1\rangle \quad .
\eeq
\nit This implies
\beq
\theta^{(n)} |m\rangle= ([m+n]_q! \, /([m]_q! \, [n]_q!) \; |m+n\rangle 
\eeq
\nit is valid for generic $q$.
Setting $m=rn+p$ as in section four and passing to the limit for $q$
a root of 1 with the aid of lemmas from section four gives
$z|rn+p\rangle =(r+1)|(r+1)n+p\rangle$. Also $\partial_z=\der^n$ leads to
$\partial_z |rn+p\rangle =|(r-1)n+p\rangle $. Indeed we can see that the 
representation space at generic $q$ in the limit acquires a  
product structure. 
Setting $|rn+p\rangle \equiv |r\, , \, p \rangle \in
{{\cal V}}_{HO} \otimes {{\cal V}}^n $, we may view $z, \partial_z$
as $z \otimes 1, \partial_z \otimes 1$ in the product space, so that
\beq
z|r\rangle =(r+1)|r+1\rangle \quad , 
\quad \partial_z |r\rangle =|r-1\rangle \quad .
\eeq
\nit Likewise we may view $\theta$, etc., as $1 \otimes \theta$ and use
\beq 
\der =1 \otimes \partial_\theta +\partial_z \otimes \theta^{(n-1)} \quad .
\eeq
\nit to express $\der$ in terms of creation and destruction operators.
There is of course a similarity transformation involved in placing
the representations considered here 
explicitly in equivalence with those in which
\beq
a|p\rangle= \Bigl( {{q^{p/2}-q^{-p/2}} \over {q^{1/2}-q^{-1/2}}} \Bigr)^{1/2} 
 \; |p-1\rangle \quad ,
\eeq
and in which the correct adjoint properties of $a^{\dagger}$ are evident.
This is discussed in \cite{DMdAPB} .

\section*{Acknowledgements}
This paper describes research supported in part by E.P.S.R.C and P.P.A.R.C. 
(UK) and by the C.I.C.Y.T (Spain).
J.C.P.B. wishes to acknowledge an FPI grant
from the CSIC and the Spanish Ministry of Education and Science.

\end{document}